\pdfoutput=1
\documentclass[11pt]{article}
\usepackage[a4paper]{geometry}
\usepackage{amsthm}
\usepackage{amsmath}
\usepackage{amssymb}
\usepackage{amsfonts}
\usepackage{graphicx}
\usepackage{graphics}
\usepackage[bf,SL,BF]{subfigure}

 \numberwithin{equation}{section}
\def\eref#1{(\ref{#1})}

\def\P{\mathbb{P}}
\def\R{\mathbb{R}}

\def\E{\mathbb{E}}

\def\B{\mathcal{B}}

\def\<{\big\langle}
\def\>{\big\rangle}

\newtheorem{Lemma}{Lemma}[section]

\newtheorem{Theorem}{Theorem}[section]
\newtheorem{Proposition}{Proposition}[section]

\newtheorem{Condition}{Condition}[section]

\theoremstyle{remark}
\newtheorem{Remark}{Remark}[section]

\theoremstyle{definition}

\theoremstyle{definition}

\newcommand{\figbox}[1]{%
  \fbox{%
    \vbox to 1in{%
    \vfil
    \hbox to 2in{%
      \hfil
      #1%
      \hfil}%
    \vfil}}}

\begin{document}
\title{Ergodicity of Langevin Processes with Degenerate Diffusion in Momentums}

\author{Nawaf Bou-Rabee\footnote{California Institute of Technology, Applied \&
Computational Mathematics,  nawaf@acm.caltech.edu}  \and Houman
Owhadi\footnote{ California Institute of Technology, Applied \&
Computational Mathematics, Control \& Dynamical Systems, MC 217-50
Pasadena , CA 91125, USA. owhadi@caltech.edu}. }

\date{\today}

\maketitle

\begin{abstract}
This paper introduces a geometric method for proving ergodicity of
degenerate noise driven stochastic processes. The driving noise is
assumed to be an arbitrary Levy process with non-degenerate
diffusion component (but that may be applied to a single degree of
freedom of the system). The geometric conditions are the approximate
controllability of the process the fact that there exists a point in
the phase space where the interior of the image of a point via  a
secondarily randomized version of the driving noise is non void.

The paper applies the method to prove ergodicity of a sliding disk
governed by Langevin-type equations (a simple stochastic rigid body
system).  The paper shows that a key feature of this Langevin
process is that even though the diffusion and drift matrices
associated to the momentums are degenerate, the system is still at
uniform temperature.
\end{abstract}
\maketitle

\section{Introduction}
This  paper is concerned with proving ergodicity of mechanical
systems governed by Langevin-type equations driven by Levy processes
and with a singular diffusion matrix applied on the momentums.

Such systems arise, for instance, when one models stochastically
forced mechanical systems composed of rigid bodies.  In such systems
one would like to introduce a certain structure to the noise and
observe its effect on the dynamics of the system. For instance, one
would like to apply stochastic forcing to a single degree of freedom
and characterize the ergodicity of the system.  The stochastic
process associated to the dynamics of these systems is in general
only weak Feller and not strong Feller.

The paper provides a concrete weak Feller (but not strong Feller)
stochastic process to illustrate this lack of regularity.  The
example is a simple mechanical system that is randomly forced and
torqued and that preserves the Gibbs measure.   In this case one
would like to determine if this Gibbs measure is the unique,
invariant measure of the system.

A new strategy based on the introduction of the \emph{asymptotically
strong Feller} property has been introduced in \cite{MR2259251}.
This paper proposes an alternative method based on two conditions:
weak irreducibility and closure under second randomization of the
stochastic forcing (see theorem \ref{ksjhskjsh223}). Our strategy is
in substance similar to the one proposed by Meyn and Tweedie for
discrete Markov Chains in chapter 7 of \cite{MeTw96}.

Although the H\"{o}rmander condition (\cite{RogWill06} 38.16) can
also be used to obtain local regularity properties of the
semi-group, hence a local strong Feller condition and ergodic
properties. The alternative approach proposed here doesn't require
smooth vector fields or manifolds, it can directly be applied to
Levy processes and (this is our main motivation) it allows for an
explicit geometric understanding of the mechanisms supporting
ergodicity.

For related previous work we refer to \cite{MR1931266},
\cite{MR1924935}, \cite{MR2259251}, \cite{MR1838749},
\cite{MR1764365}, \cite{MR1705589} and \cite{MR1685893}.

\section{General set up.}
Let $(X_t)_{t\in \R^+}$ be a Markov stochastic process on a
(separable) manifold $M$ with model space $\mathbb{R}^n$.

Let $(\omega_t)_{0\leq t}$ be $p$-dimensional Levy process, i.e. a
  stochastic process on $\R^p$ that has has independent increments, is stationary, is stochastically
 continuous and such that (almost surely) trajectories are continuous from the left and with limits
from the right.

We assume that there exists a family deterministic mappings (indexed
by $0 \leq t$) $F_{t}: M \times ([0,t]\rightarrow \mathbb{R}^p) \to
M$ such that
\begin{equation}
\begin{cases}
X_t= F_{t-s}\big(X_s, (\omega_{s'}-\omega_s)_{s \leq s' \leq t}
\big)
\end{cases}
\end{equation}

Recall that the first three condition defining a Levy process mean
that $(\omega_t-\omega_s)_{t\geq s}$ is independent of
$(\omega_{s'})_{0\leq s'\leq s})$, the law of $\omega_t-\omega_s$
depends only on $t-s$ and $\lim_{s\rightarrow
0}\P[|\omega_{s+t}-\omega_t|\geq \epsilon]= 0$.

Recall also \cite{RogWill06, Win08} that since $\omega$ is a Levy
process, there exists a $\gamma \in \R$, a constant $p\times p$
matrix $\sigma$, a standard $p$-dimensional Brownian Motion
$(B_t)_{t\geq 0}$ and $(\Delta_t)_{t\geq 0}$ an independent Poisson
process of jumps with intensity of measure $dt \times \nu(dx)$ on
$dt \times \R^p$ (such that
$\int_{\R^p}\min(1,|z^p|)\nu(dz)<\infty$) such that
\begin{equation}\label{hjehej3errr}
\omega_t=\gamma t + \sigma B_t + C_t+M_t
\end{equation}
Where $C_t=\sum_{s\leq t} \Delta_s 1_{|\Delta_s|> 1}$ is a compound
Poisson point process (of jumps of norm larger than one) and
\begin{equation}
M_t=\lim_{\epsilon \downarrow 0}\Big(\Delta_s 1_{\epsilon<|\Delta
s|\leq 1}- t \int_{z\in \R^p\,:\, \epsilon<|z|\leq 1}z \nu(dz)\Big)
\end{equation}
is a martingale (of small jumps compensated by a linear drift).
Recall also that any process that can be represented as
\eref{hjehej3errr} is a $p$-dimensional Levy process, in particular
a $p$-dimensional Brownian Motion is a Levy process. In this paper,
the only assumption on the stochastic forcing $\omega$ will be the
following one:

\begin{Condition}\label{cond_omega}
$\sigma$ is non degenerate (has a non null determinant).
\end{Condition}

We will then prove the ergodicity of $X_t$ based on the following
geometric conditions on $F$.

\begin{Condition}\label{shsksjsa}
$X_t$ is  approximately  controllable, i.e., for all $A,B\in M$ and
$\epsilon >0$ there exists $t>0$ and $\phi \in C^0([0,t], \R^p)$
 so that $F_{t}\big(A,(\phi_s-\phi_0)_{0\leq s \leq t} \big) \in
\mathcal{B}(B,\epsilon)$.
\end{Condition}

\begin{figure}[htbp]
\begin{center}
\includegraphics[scale=0.5,angle=0]{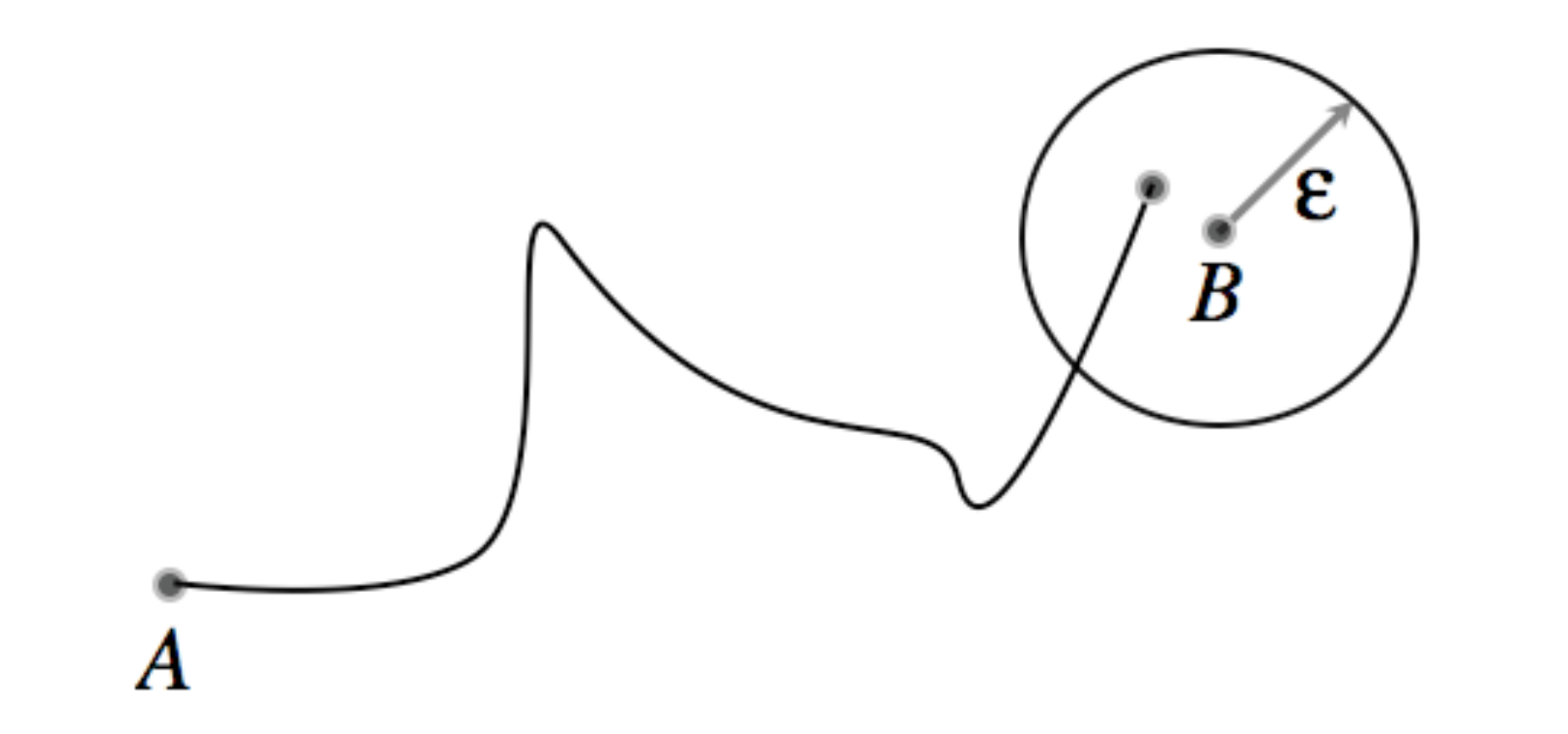}
\caption{\footnotesize  {\bf Approximate Controllability Condition.}
The condition states that given $A, B \in M$ and $\epsilon>0$, there
exists $t>0$ and $\phi \in C^0([0,t], \R^p)$
 so that $F_{t}\big(A,(\phi_s-\phi_0)_{0\leq s \leq t} \big) \in
\mathcal{B}(B,\epsilon)$.
 } \label{fig:weakirreducibility}
\end{center}
\end{figure}

This condition is illustrated in Fig.~\ref{fig:weakirreducibility}.

\begin{Condition}\label{shsksjsaN}
For all  $0 \leq t$, the mapping $(x,\phi) \mapsto F_{t}\big(x,
(\phi_{s}-\phi_0)_{0 \leq s \leq t} \big)$ is continuous with
respect to the norm $\|x-y\|+\|\phi-\psi\|$ where
$\|\phi-\psi\|:=\sup_{0 \leq s \leq
t}|\phi_s-\phi_0-(\psi_s-\psi_0)|$.
\end{Condition}

Let $\phi,\varphi^1,\ldots,\varphi^n$ be $n+1$ deterministic
continuous mappings from $[0,t]$ onto $\R^n$ equal to $0$ at time
$0$. For $\lambda \in \R^n$, write
 \begin{equation}
G(a,\phi,\lambda):=F_t\big(a,(\phi_s+\sum_{i=1}^n \lambda_i
\varphi^i_s)_{0\leq s \leq t}\big)
\end{equation}

\begin{Condition}\label{shssesswedeksjs}
There exists $x_0\in M$ and $t>0$, such that in a neighborhood of
$(x_0,0,0)$:

\begin{itemize}
\item $(x,(\phi)_{0\leq s \leq t},\lambda)\rightarrow
G(x,(\phi)_{0\leq s \leq t},\lambda)$ is
 differentiable in $\lambda$.
 \item  $\nabla_\lambda G$ is invertible and uniformly bounded.
 \item  $(\nabla_\lambda G)^{-1}$
is uniformly bounded.
 \end{itemize}
\end{Condition}

\section{Main theorem}
\begin{Theorem}\label{ksjhskjsh223}
Consider a stochastic process $X_t$ on a manifold $M$ that satisfies
conditions \ref{shsksjsa}, \ref{shsksjsaN}, \ref{shssesswedeksjs},
\ref{cond_omega} and admits an invariant measure $\mu$. Let $P_t$ be
the semigroup associated to $X$.  Then
\begin{itemize}
\item $\mu$ is ergodic and weakly mixing with respect to $P_t$.
\item $\mu$ is the unique invariant measure of $X$.
\end{itemize}
\end{Theorem}

\begin{proof}
We will need the following two lemmas on the Levy process $\omega$.

\begin{Lemma}\label{shsksjdzs}
Assume that $\omega$ satisfies condition \ref{cond_omega}. Let
$0\leq s <t$ and $\phi \in C^0([0,t], \R^p)$ be arbitrary. The laws
of $\big(\omega_{s}-\omega_{0} \big)_{0\leq s \leq t}$ and
$\big(\omega_{s}-\omega_{0}-(\phi_{s}-\phi_0) \big)_{0\leq s \leq
t}$ are absolutely continuous with respect to each other.
\end{Lemma}
\begin{proof}
Lemma \ref{shsksjdzs} follows by applying Girsanov's theorem to the
diffusive component ($B$) of $\omega$.
\end{proof}

\begin{Lemma}\label{shsksjsb}
Assume that $\omega$ satisfies condition \ref{cond_omega}. Let $\phi
\in C^0([0,t], \R^p)$. For all $\epsilon>0$, the inequality
$\P\big[\sup_{0\leq s \leq
t}|\phi_{s}-\phi_0-(\omega_{s}-\omega_0)|<\epsilon\big]>0$ holds
almost surely.
\end{Lemma}
\begin{proof}
Let $\epsilon>0$. Let $(\gamma,\sigma,\nu)$ be the Levy-Khintchine
characteristics of $\omega$. Let $\eta>0$ such that
\begin{equation}
\int_{z\in \R^p\,:\, 0<|z|\leq \eta}z^2\nu(dz)<\frac{\epsilon^4}{16}
\end{equation}
Observe that \cite{Win08} can be written as
\begin{equation}\label{hjehej3ejjrrr}
\omega_t=\gamma^\eta t + \sigma B_t + C_t^\eta+M_t^\eta
\end{equation}
where
\begin{equation}\label{hjedkdhej3errr}
\gamma^\eta=\gamma-\int_{z\in \R^p\,:\, \eta <|z|\leq 1}z \nu(dz),
\end{equation}
and
\begin{equation}\label{hjedkdhghej3errr}
C_t^\eta=\sum_{s\leq t} \Delta_s 1_{|\Delta_s|> \eta}
\end{equation}
 is a
compound Poisson point process (of jumps of norm larger than one)
and
\begin{equation}
M_t^\eta=\lim_{\epsilon \downarrow 0}\Big(\Delta_s
1_{\epsilon<|\Delta s|\leq \eta}- t \int_{z\in \R^p\,:\,
\epsilon<|z|\leq \eta}z \nu(dz)\Big)
\end{equation}
is a martingale (of small jumps compensated by a linear drift).
Observe that with strictly positive probability $\exp(-t
\nu(|z|>\eta))$, $C_t^\eta$ is uniformly equal to 0 over $[0,t]$.
Furthermore by the Martingale maximal inequality
\begin{equation}
\E[\sup\{(M^\eta_s)^2\,:\, 0 \leq s \leq t\}] \leq 4
\E[(M^\eta_t)^2]
\end{equation}
and using (\cite{Win08})
\begin{equation}
\E[(M^\eta_t)^2]=\int_{z\in \R^p\,:\, 0<|z|\leq \eta}z^2\nu(dz)
\end{equation}
and Chebyshev's inequality we obtain that
\begin{equation}
\P[\sup_{0\leq s\leq t}|M_s|\geq \frac{\epsilon}{2}] \leq
\frac{\epsilon}{2}
\end{equation}
hence
\begin{equation}
\P[\sup_{0\leq s\leq t}|M_s|< \frac{\epsilon}{2}] \geq
1-\frac{\epsilon}{2}.
\end{equation}
We conclude the proof of lemma \ref{shsksjsb} by applying Schilder's
theorem  to  $B_t$ and using the fact that $\sigma$ is not
degenerate.
\end{proof}

Let us now prove that $\mu$ is ergodic. Let $A\in \B(M)$ be an
invariant set of positive $\mu$-measure, i.e.,
\begin{equation}
P_t 1_{A}=1_{A},\quad\text{for every}\quad t\geq 0,\quad \mu-a.s.
\end{equation}
and $\mu(A)>0$. We will prove that $\mu(A)=1$.  Assume $0<\mu(A)<1$.
Then $A^c$, which is also an invariant set, has strictly positive
measure, i.e. $\mu(A^c)>0$. Now let us prove the following lemma

\begin{Lemma}\label{ssnkje8i33}
If $0<\mu(A)<1$ then
\begin{itemize}
\item For all $y\in M$ and $\epsilon>0$, $\mu(A \cap
\mathcal{B}(y,\epsilon))>0$.
\item For all $y\in M$ and $\epsilon>0$, $\mu(A^c \cap
\mathcal{B}(y,\epsilon))>0$.
\end{itemize}

\end{Lemma}
\begin{proof}
 We will restrict the proof to $A$. Since
$\mu(A)>0$ there exists $x_0>0$ such that for all $\epsilon>0$,
$\mu(A\cap \mathcal{B}(x_0,\epsilon))>0$ (otherwise one would get
$\mu(A)=0$ by covering the separable manifold $M$ with a countable
number of balls such that $\mu(A\cap \mathcal{B}(x,\epsilon_x))=0$).
Assume that there exists $y_0\in M$ and $\epsilon>0$ such that
$\mu(A \cap \mathcal{B}(y_0,\epsilon))=0$.  Since $X_t$ is weakly
controllable (condition \ref{shsksjsa}) there exists $t>0$ and $\phi
\in C^0([0,t], \R^p)$
 so that $F_{t}\big(x_0,(\phi_s-\phi_0)_{0\leq s \leq t} \big) \in
\mathcal{B}(y_0,\frac{\epsilon}{2})$. From the continuity condition
\ref{shsksjsaN} on $F$ and the Schilder type lemma \ref{shsksjsb}
imply that there exists $\epsilon'>0$ such that
\begin{equation}\label{eqgjdeid}
 \text{for all}\quad x\in \mathcal{B}(x_0,\epsilon'),\quad
\P\big[ F_{t}\big(x,(\phi_s-\phi_0)_{0\leq s \leq t} \big) \in
\mathcal{B}(y_0,\epsilon)\big]>0. \end{equation}
 Write $P_t$ the
semi-group associated with $X_t$. Equation \eref{eqgjdeid} leads to
a contradiction with the fact that
\begin{equation}
\int_{A} P_t(x,A) \mu(dx)=\mu(A).
\end{equation}
since $\mu\big(A\cap \mathcal{B}(x_0,\epsilon') \big)>0$ and for all
$x\in \mathcal{B}(x_0,\epsilon')$, $P_t(x,A)<1$.

\end{proof}
From condition \ref{shssesswedeksjs} there exists $x_0\in M$ and
$t,\epsilon, \alpha, \delta, K>0$ and such that for $x\in
B(x_0,\epsilon)$, $\|\phi\|_{L^\infty(0,t)}<\alpha$ and $\lambda \in
(-\delta,\delta)^n$, $G(x,(\phi)_{0\leq s \leq t},\lambda)$ is
differentiable in $\lambda$, $|\nabla_{\lambda} G|\leq K$ and
$\big|(\nabla G)^{-1}\big|\leq K$. It follows from the condition
\ref{shssesswedeksjs} and the continuity condition \ref{shsksjsaN}
that $\epsilon'\in (0,\epsilon)$ can be chosen small enough so that
there exists $z\in M$, $0<\alpha'<\alpha$, $0<\epsilon_z$ such that
for $\|\phi\|_{L^\infty(0,t)}<\alpha'$ we have for all $a,b \in
B(x_0,\epsilon')$,
\begin{equation}\label{ksljkehejh3e}
\mathcal{B}(z,\epsilon_z) \subset G(a,(\phi)_{0\leq s \leq
t},(-\delta,\delta)^n)\cap G(b,(\phi)_{0\leq s \leq
t},(-\delta,\delta)^n).
\end{equation}
Equation \eref{ksljkehejh3e} is illustrated in
Fig.~\ref{fig:closure}.
\begin{figure}[htbp]
\begin{center}
\includegraphics[scale=0.5,angle=0]{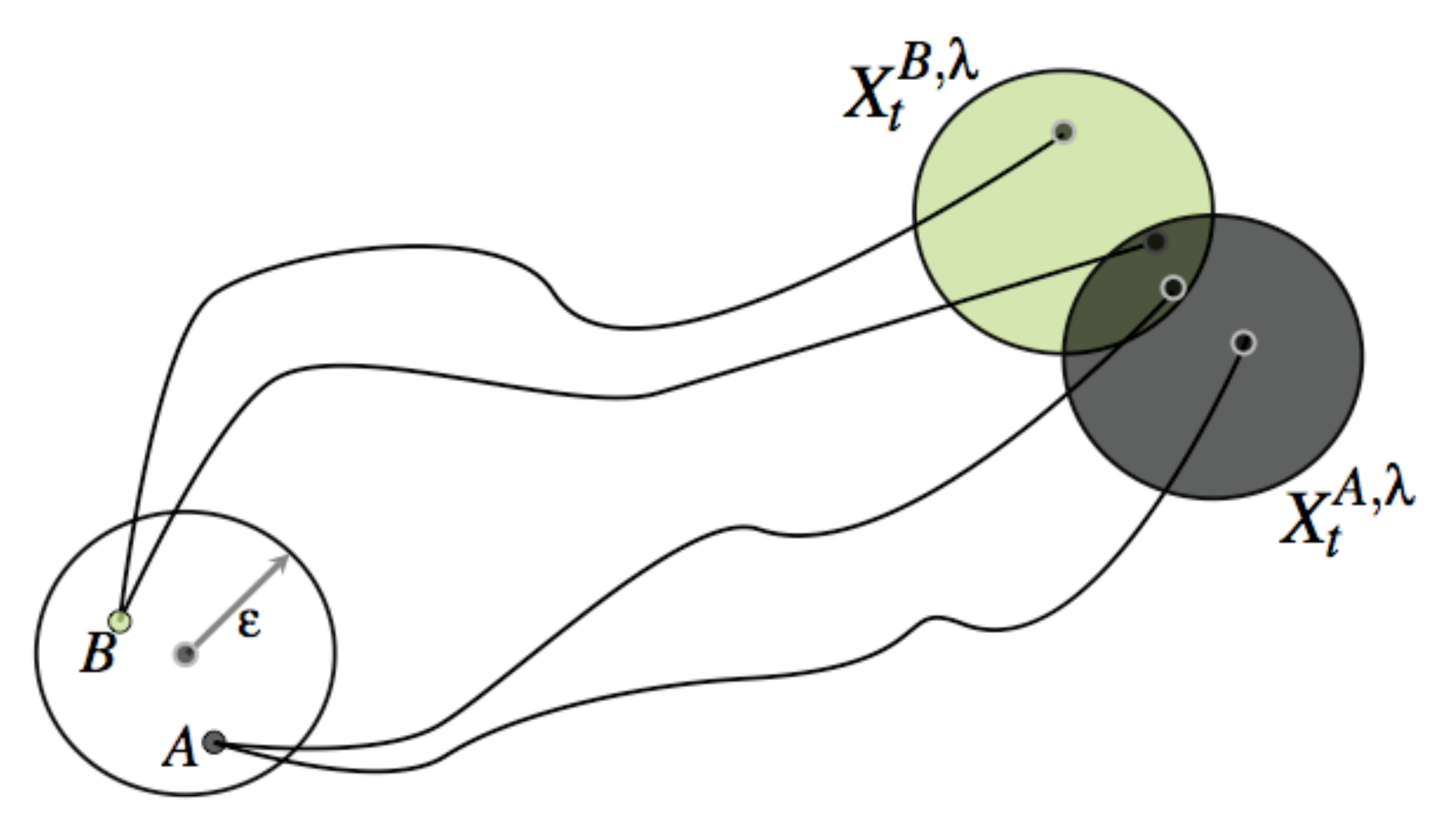}
\caption{\footnotesize  {\bf Closure Under Second Randomization
Condition Illustrated.} This condition states that under a second
randomization of the noise via $\lambda$, the interior of the
intersection of the range  of $G(A, \lambda)$ (image of
$(-\delta,\delta)^n$ by $\lambda \rightarrow G(A,.)$) and $G(B,
\lambda)$ is not void.
 } \label{fig:closure}
\end{center}
\end{figure}

 Let $T>t$.
 From the previous lemma there exists $a \in
\mathcal{B}(x_0,\epsilon') \cap A$ and $b \in
\mathcal{B}(x_0,\epsilon') \cap A^c$ such that $P_T(a,A)=1$ and
$P_T(b,A^c)=1$.  Set $X_t^a$ ($X_t^b$) to be the process $X_t$
started from the point $a \in M$ ($b\in M$) and set $\P_a$ to be the
measure of probability associated to $X_t^a$. We obtain from the
Markov property that
\begin{equation}
\E\big[P_{T-t}(X^a_t,A)\big]=1 \quad \text{and}\quad
\E\big[P_{T-t}(X^b_t,A)\big]=0 \text{.}
\end{equation}
Write
\begin{equation}
X^{a,\lambda}:=F_t\big(a,(\omega_s-\omega_0+\sum_{i=1}^n \lambda_i
\varphi^i_s)_{0\leq s \leq t}\big)
\end{equation}
The Girsanov type lemma \ref{shsksjdzs} implies that the laws of
$X^a$ and $X^{a,\lambda}$ are absolutely continuous with respect to
each other. Hence for all $\lambda \in (-\delta,\delta)^n$,
\begin{equation}
\E\big[P_{T-t}(X^{a,\lambda}_t,A)\big] =1 \quad \text{and}\quad
\E\big[P_{T-t}(X^{b,\lambda}_t,A)\big]=0.
\end{equation}
Which leads to
\begin{equation}\label{ksjhslkh3s}
\delta^{-2n}\int_{[-\delta,\delta]^n}\E\big[P_{T-t}(X^{a,\lambda}_t,A)\big]\,d\lambda
=1 \quad \text{and}\quad
\delta^{-2n}\int_{[-\delta,\delta]^n}\E\big[P_{T-t}(X^{b,\lambda}_t,A)\big]\,d\lambda=0.
\end{equation}
Let $\Omega_I$ be the event $\|\omega\|_{L^\infty(0,t)}<\alpha'$.
Observe that from the Schilder type lemma \ref{shsksjsb} the measure
of probability of
 $\Omega_I$ is  strictly positive. It follows from \eref{ksjhslkh3s}
 and \eref{ksljkehejh3e}
 that
\begin{equation}
\delta^{-2n}\int_{[-\delta,\delta]^n}\E\big[1_{\Omega_I}1_{X^{a,\lambda}_t\in
\mathcal{B}(z,\epsilon_z) }P_{T-t}(X^{a,\lambda}_t,A)\big]\,d\lambda
>0
\end{equation}
Using the change of variable $y=X^{a,\lambda}_t$ we obtain from
\eref{ksljkehejh3e} that
\begin{equation}
\E\big[1_{\Omega_I} \int_{\mathcal{B}(z,\epsilon_z)}
P_{T-t}(y,A)\frac{dy}{|\nabla_\lambda X^{a,\lambda}|\circ
(X^{a,\lambda})^{-1}(y)}\big]
>0
\end{equation}
Hence
\begin{equation}
\E\big[1_{\Omega_I} \int_{\mathcal{B}(z,\epsilon_z)}
P_{T-t}(y,A)\frac{|\nabla_\lambda X^{b,\lambda}|\circ
(X^{b,\lambda})^{-1}(y)}{|\nabla_\lambda X^{a,\lambda}|\circ
(X^{a,\lambda})^{-1}(y)} \frac{dy}{|\nabla_\lambda
X^{b,\lambda}|\circ (X^{b,\lambda})^{-1}(y)}\big]
>0
\end{equation}
We deduce from equation \eref{ksljkehejh3e} and the fact that
$\frac{|\nabla_\lambda X^{b,\lambda}|\circ
(X^{b,\lambda})^{-1}(y)}{|\nabla_\lambda X^{a,\lambda}|\circ
(X^{a,\lambda})^{-1}(y)} $ is bounded from below by $K^{-2}$ that
\begin{equation}
\E\big[1_{\Omega_I} \int_{\mathcal{B}(z,\epsilon_z)} P_{T-t}(y,A)
\frac{dy}{|\nabla_\lambda X^{b,\lambda}|\circ
(X^{b,\lambda})^{-1}(y)}\big]
>0.
\end{equation}
However a similar computation leads from  \eref{ksjhslkh3s}
 and \eref{ksljkehejh3e} to
 \begin{equation}
\E\big[1_{\Omega_I} \int_{\mathcal{B}(z,\epsilon_z)} P_{T-t}(y,A)
\frac{dy}{|\nabla_\lambda X^{b,\lambda}|\circ
(X^{b,\lambda})^{-1}(y)}\big] =0.
\end{equation}
Hence a contradiction. Thus $\mu$ must be ergodic. Let us now prove
that $\mu$ is the unique invariant measure. Assume that
$\mu'\not=\mu$ is also invariant with respect to the semigroup
$P_t$. By the argument presented above $\mu'$ is ergodic and it
follows from Proposition 3.2.5 of \cite{PrZa96} that $\mu$ and
$\mu'$ are singular and it is easy to check from the argument
presented above that this can't be the case (the proof is similar to
the one given in theorem 4.2.1 of \cite{PrZa96}). Hence $\mu$ is the
unique invariant distribution. The proof of the fact that $\mu$ is
weakly mixing follows from theorem 3.4.1 of \cite{PrZa96} and is
similar to the one given at page 44 of \cite{PrZa96} (theorem
4.2.1).
\end{proof}

\section{Sliding Disk at Uniform Temperature.}
\begin{figure}[htbp]
\begin{center}
\includegraphics[scale=0.5,angle=0]{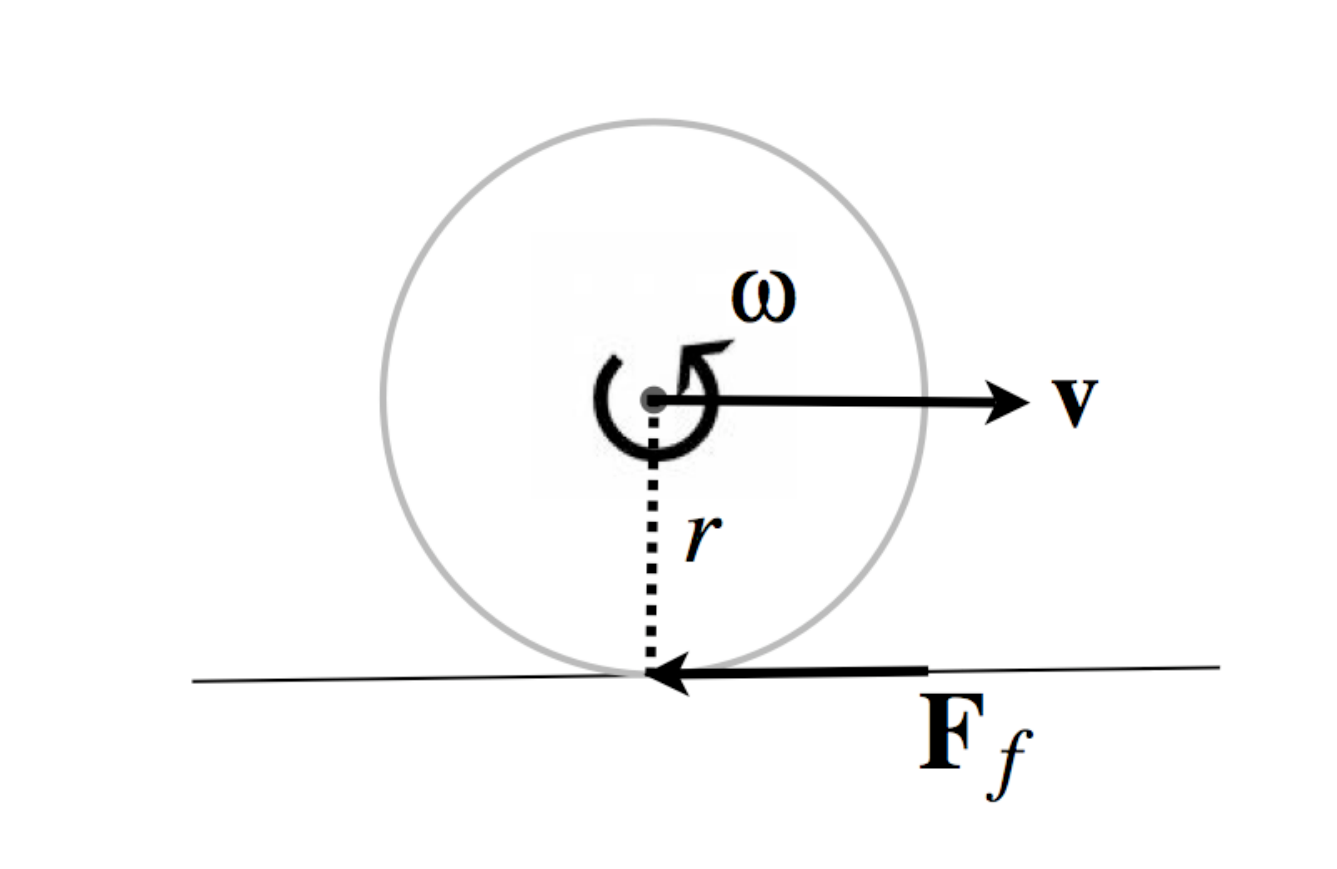}
\caption{\footnotesize  {\bf Sliding Disk.}   Consider a sliding
disk of radius $r$ that is free to translate and rotate on a
surface.   We assume the disk is in sliding frictional contact with
the surface.  The configuration space of the system is
$\operatorname{SE}(2)$, but with the surface constraint the
configuration space is just $\mathbb{R} \times
\operatorname{SO}(2)$. } \label{fig:slidingdisk}
\end{center}
\end{figure}

Consider a disk on a surface as shown in Fig.~\ref{fig:slidingdisk}
\cite{BoOw2007c}. The disk is free to slide and rotate.  We assume
that one rescales position its radius and time by some
characteristic frequency of rotation or other time-scale. The
dimensionless Lagrangian is given by
\begin{equation}
L(x, v, \theta, \omega) = \frac{1}{2} v^2 + \frac{\sigma}{2}
\omega^2 - U(x)
\end{equation}
where $v$ stands for the velocity of the center of mass, $\omega$
the angular velocity of the disk and $\sigma$ is a strictly positive
dimensionless constant given by $\sigma := J/(m r^2)$ (where $r$ is
the radius of the disk, $m$ is its mass and $J$ its moment of
inertia). $U: \mathbb{R} \to \mathbb{R}$ is an arbitrary periodic
potential which is assumed to be smooth, and of period one.

\begin{figure}[htbp]
\begin{center}
\includegraphics[scale=0.4,angle=0]{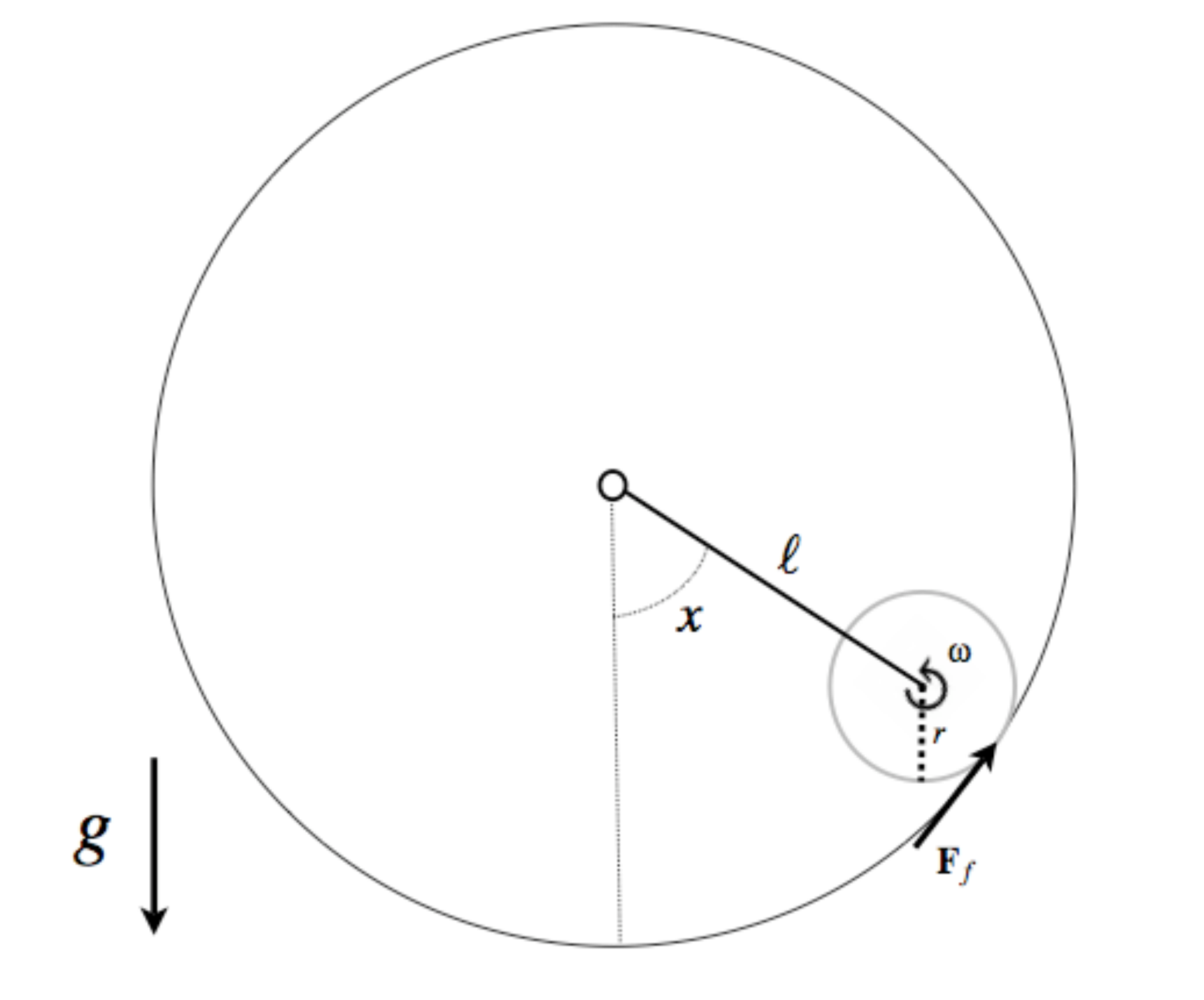}
\caption{\footnotesize  {\bf Ballistic Pendulum.}   If the
dimensionless potential is $U =  \cos(x)$, then the sliding disk is
simply a pendulum in which the bob in the pendulum is replaced by a
disk and the pendulum is placed within a cylinder as shown.     }
\label{fig:ballisticpendulum}
\end{center}
\end{figure}

The contact with the surface is modeled using a sliding friction
law.  For this purpose we introduce a symmetric matrix $\mathbf{C}$
defined as,
\[
\mathbf{C} = \begin{bmatrix} 1 &  1/\sigma \\
                             1/\sigma & 1/\sigma^2  \end{bmatrix} \text{.}
\]
Observe that $\mathbf{C}$ is degenerate since the frictional force
is actually applied to only a single degree of freedom, and hence,
one of its eigenvalues is zero.  In addition to friction a white
noise is applied to the same degree of freedom to which friction is
applied. The governing stochastic differential equations are

\begin{equation}\label{sjhgsh3}
\begin{cases}
 d x &= v dt  \\
 d \theta &= \omega dt  \\
\begin{bmatrix} d v \\   d  \omega \end{bmatrix} &=
\begin{bmatrix} - \partial_x U  \\ 0 \end{bmatrix} dt
- c \mathbf{C} \begin{bmatrix}   v \\    \sigma \omega
\end{bmatrix} dt + \alpha \mathbf{C}^{1/2} \begin{bmatrix} d B_v \\
d B_{\omega}
\end{bmatrix} \text{.}
\end{cases}
 \end{equation}
 where $\mathbf{C}^{1/2}$ is the matrix square root of $\mathbf{C}$.    The matrix square root is easily computed by diagonalizing $\mathbf{C}$ and computing square roots of the diagonal entries (eigenvalues of $\mathbf{C}$) as shown:
 \[
\mathbf{C}^{1/2}
 = \frac{\sigma}{\sqrt{\sigma^2+1}} \mathbf{C} \text{.}
 \]
Write $X:=(x,\theta,v,\omega)$. It easy to check that the  Gibbs
distribution
\begin{equation}
\mu(d\xi):=\frac{e^{-\beta E}}{Z} dX
\end{equation}
is invariant for \eref{sjhgsh3}, where  $\beta = 2 c / \alpha^2$,
$Z:=\int e^{-\beta E}dX$, and $E$ is the energy of the mechanical
system and is given by
\[
E := \frac{1}{2} v^2 + \frac{1}{2} \sigma \omega^2  + U(x) \text{.}
\]
Define
\begin{equation}
Y:=\begin{pmatrix}-x+\sigma \theta \\ x+\theta \end{pmatrix}
\end{equation}
The system \eref{sjhgsh3} can be written
\begin{equation}\label{kjjsksssesshs}
\begin{cases}
\dot{Y_1}(t)=\dot{Y_1}(0)+\int_0^t \partial_x U(\frac{\sigma
Y_2-Y_1}{\sigma+1})\,ds \\
\dot{Y_2}(t)=\dot{Y_2}(0)-\int_0^t \partial_x U(\frac{\sigma
Y_2-Y_1}{\sigma+1})\,ds-c\gamma (Y_2(t)-Y_2(0))+\bar{\alpha}
\sqrt{2}B_t
\end{cases}
\end{equation}
where $\gamma=(\sigma+1)/\sigma$,
$\bar{\alpha}=\alpha(\sigma+1)/\sqrt{\sigma^2+1}$ and
$B:=(B_v+B_\omega)/\sqrt{2}$ is a one dimensional Brownian Motion.
Observe that condition \ref{cond_omega} is satisfied with
$\omega=B$, $p=1$ and $\sigma=(1)$.

 Observe also that if $U$ is a constant then the quantity
$-v+\sigma \omega$ is conserved and the system \eref{sjhgsh3} can't
be ergodic. Let us assume that $U$ is not constant, our purpose is
to prove that the Gibbs distribution $\mu$ is ergodic with respect
to the stochastic process $X$.

\begin{Remark}
Observe that when $U$ is not constant over a non void open subset of
$\R$ (say $(-\frac{1}{4},\frac{1}{4})$), $Y$ needs to travel a
distance that is uniformly (in $\epsilon$) bounded from below by a
strictly positive amount to get from
$(Y_1,Y_2,\dot{Y_1},\dot{Y_2})=(0,0,0,0)$ to the domain
$\dot{Y_1}>\epsilon$. It follows that in that situation that the
process $Y$ and hence $X$ is not strong Feller and theorems
requiring this property can't be applied.
\end{Remark}

\begin{Remark}
Observe also that  the condition $\partial_{x}^2 U\not=0$ in a
neighborhood of $x_0$ doesn't guarantee that $Y$ is strongly Feller
in that neighborhood. For instance observe that $\partial_{x}^2
U(x_0)\not=0$ and $\partial U(x_0)>0$ imply that the drift on $Y_1$
is uniformly bounded by a strictly positive constant on a
neighborhood of $(0,\frac{\sigma+1}{\sigma} x_0)$ it follows that
$\P_\epsilon {(y_1,y_2)}[Y_1<0]$ is discontinuous in the
neighborhood of  $(0,\frac{\sigma+1}{\sigma} x_0)$ ($\epsilon$)
close to the line $y_1=0$.
\end{Remark}

We believe that the system $Y$ is asymptotically strong Feller so
one could in principle obtain the ergodicity of $\mu$ by controlling
the semi-group associated to $Y$ as it is suggested in
\cite{MR2259251}. We propose an alternative method based on the
controllability of the ODE associated to $Y$ and theorem
\ref{ksjhskjsh223}. We believe that it is much simpler to control
the geometric properties of the ODE associated to $X$ rather than
the gradient of its semigroup.

One can also check that the generator of $Y$ satisfies a local
H\"{o}rmander condition (\cite{RogWill06} 38.16) at a point $x_0$
such that $\partial_{x}^2 U(x_0)\not=0$ so an alternative method to
prove ergodicity would be to use that condition to obtain a local
regularity of the semi group associated to $U$. Here we propose an
alternative method which doesn't require $U$ to be smooth and which
can be applied with Levy processes.

\begin{Theorem} \label{thm:slidingdiskerggodic}
Assume that $U$ is not constant. Then the Gibbs measure $\mu$ is
ergodic and strongly mixing with respect to the stochastic process
$X$ \eref{sjhgsh3}. Furthermore, it is the unique invariant
distribution of $X$.
\end{Theorem}

First let us prove that codition \ref{shsksjsa} is satisfied by $X$.
\begin{Lemma}
Assume $U$ is not constant. Then $Y$ is approximately controllable.
\end{Lemma}
\begin{proof}
Since $U$ is not constant, there exists $t_1>0$ such that for
$t_i\geq t_1$ there exists a smooth path $Y$ such that
$Y_1(0)=-x_1+\sigma \theta_1$, $Y_2(0)=x_1+\theta_1$,
$\dot{Y_1}(0)=-v_1+\sigma \omega_1$, $\dot{Y_2}(0)=v_1+\omega_1$,
$Y_1(t_i)=-x_2+\sigma \theta_2$, $\dot{Y_1}(t_i)=-v_2+\sigma
\omega_2$ and
\begin{equation}\label{kjsejskshs}
\frac{d^2Y_1}{dt^2}=\partial_x U \left(\frac{\sigma
Y_2-Y_1}{\sigma+1} \right)
\end{equation}
Take $t_2:=t_i+\frac{\min(\epsilon,1)}{10 (\|\partial_x
U\|_{L^\infty}+1+|\frac{d}{dt}Y_1(t_i)|)}$ and interpolate smoothly
$Y_2$ between $Y_2(t_i)$ (obtained from the control problem
\eref{kjsejskshs}) and
\begin{equation}
\begin{pmatrix}Y_2(t_2) \\ \frac{dY_2}{dt}(t_2) \end{pmatrix}= \begin{pmatrix}x_2+\theta_2 \\
v_2+\omega_2 \end{pmatrix}
\end{equation}
Observe that the extension of $Y_1$ to $(t_i,t_2]$ as a solution of
\eref{kjsejskshs} satisfies
\begin{equation}
\Big|\begin{pmatrix}Y_1(t_2)-Y_1(t_i)\\
\frac{dY_1}{dt}(t_2)-\frac{dY_1}{dt}(t_i) \end{pmatrix}\Big|\leq
\frac{\epsilon}{5}
\end{equation}
Taking $\phi$ be the smooth curve defined by $\phi(0)=0$ and
\begin{equation}\label{kjjsdddskshs}
\frac{d^2Z_2}{dt^2}=-\partial_x U \left(\frac{\sigma
Z_2-Z_1}{\sigma+1} \right)-c\gamma \frac{dZ_2}{dt}+\bar{\alpha}
\sqrt{2}\frac{d \phi}{dt}
\end{equation}
completes the proof.
\end{proof}

\begin{proof}
The proof that $X$ satisfies condition \ref{shsksjsaN} is a standard
application of  Gronwall's  lemma. Observe that  the semi-group
associated to $X$ is not strongly irreducible and never equivalent
to $\mu$ because $|(-v+\sigma \omega)(t)-(-v+\sigma \omega)(0)|\leq
\|\partial_x U\|_{L^\infty} t$.  Let us now show that condition
\ref{shssesswedeksjs} is satisfied.

 Write $\xi$ the stochastic process defined by
\begin{equation}\label{kjjsksssessgouigyhs}
\begin{cases}
\dot{\xi_1}(t)=\dot{\xi_1}(0)+\int_0^t \partial_x U(\frac{\sigma
\xi_2-\xi_1}{\sigma+1})\,ds \\
\dot{\xi_2}(t)=\dot{\xi_2}(0)+\bar{\alpha} \sqrt{2}B_t
\end{cases}
\end{equation}
To prove that $Y$ satisfies condition \ref{shssesswedeksjs} it is
sufficient to show that $\xi$ satisfies condition
\ref{shssesswedeksjs}.

Since $U$ is smooth and not constant, there exists a point $x^0\in
[0,1)$, $\epsilon,C>0$ such that for $x\in B(x^0,\epsilon)$,
$\partial_x^2 U>C$. Let $\zeta$ be a point of the phase space such
that $\frac{\sigma \zeta_2-\zeta_1}{\sigma+1}=x^0$ and
$\dot{\zeta}_1=\dot{\zeta}_2=0$. Let $0<\epsilon'<\epsilon/100$ and
$a\in B(\zeta,\epsilon')$.

Let $\varphi_1,\ldots,\varphi_4$ be $4$ continuous mappings from
$\R^+$
 onto $\R$, equal to zero at time zero. For $\lambda\in \R^4$ we
write $\xi^\lambda$ the solution of
\begin{equation}\label{kjjsksshghghfsesshs}
\begin{cases}
\dot{\xi_1^\lambda}(t)=\dot{a}_1+\int_0^t \partial_x U(\frac{\sigma
\xi_2^\lambda-\xi_1^\lambda}{\sigma+1})\,ds \\
\dot{\xi_2^\lambda}(t)=\dot{a}_2+\bar{\alpha}
\sqrt{2}\sum_{i=1}^4 \lambda_i \varphi_i(t)\\
\xi_1^\lambda(t)=a_1 +\int_0^t \dot{\xi_1^\lambda}(s)\,ds\\
\xi_2^\lambda(t)=a_2 +\int_0^t \dot{\xi_2^\lambda}(s)\,ds
\end{cases}
\end{equation}
It follows that
\begin{equation}\label{kjjsddssesshs}
\begin{cases}
\dot{\xi_1^\lambda}(t)-\dot{\xi_1^0}(t)=  \int_0^t (\frac{\sigma
\xi_2^\lambda-\xi_1^\lambda}{\sigma+1}-\frac{\sigma
\xi_2^0-\xi_1^0}{\sigma+1}) \int_0^1 \partial_x^2 U\big(\frac{\sigma
\xi_2^0-\xi_1^0}{\sigma+1}+\alpha (\frac{\sigma
\xi_2^\lambda-\xi_1^\lambda}{\sigma+1}-\frac{\sigma
\xi_2^0-\xi_1^0}{\sigma+1}) \big)(s)\,ds\,d\alpha \\
\dot{\xi_2^\lambda}(t)-\dot{\xi_2^0}(t)=\sum_{i=1}^4 \lambda_i
 \big(\bar{\alpha}
\sqrt{2}\varphi_i(t)\big)\\
\xi_1^\lambda(t)-\xi_1^0(t)=\int_0^t (\dot{\xi_1^\lambda}(s)-\dot{\xi_1^0}(s))\,ds\\
\xi_2^\lambda(t)-\xi_2^0(t)=\sum_{i=1}^4 \lambda_i \big(\bar{\alpha}
\sqrt{2} \int_0^t\varphi_i(s)\,ds\big)
\end{cases}
\end{equation}
Writing $\eta$ the solution of
\begin{equation}\label{kjjsddsssdesshs}
\begin{cases}
\dot{\eta}(t)+ \frac{\partial_x^2 U\big(x^0\big)}{\sigma+1}\int_0^t
\eta(s)\,ds= \partial_x^2 U\big(x^0\big) \frac{\sigma }{\sigma+1}
\int_0^t (\xi_2^\lambda-\xi_2^0)(s)
\,ds \\
\eta=\int_0^t \dot{\eta}(s)\,ds
\end{cases}
\end{equation}
we obtain that up to the first order in $\lambda$, and at the order
$0$ in $\epsilon'$ and $t$,
$\big(\dot{\xi_1^\lambda}(t)-\dot{\xi_1^0}(t),
\xi_1^\lambda(t)-\xi_1^0(t)\big) $ can be approximated by
$(\dot{\eta}(t),\eta(t))$. It follows that $\xi^\lambda_t-\xi^0_t$
can be written as $M(\lambda,t) \lambda$ where $M(\lambda,t)$ is
continuous in $t$ and $\lambda$ in the neighborhood of $0$.
Moreover, $\varphi_1,\ldots,\varphi_4$ can be chosen so that $M$,
and $M^{-1}$ are uniformly bounded in that neighborhood. Choosing
  $0<  \delta \ll 1$ and  $0<\epsilon' \ll  \delta t \ll 1$ implies
condition  \ref{shssesswedeksjs}. By invoking theorem
\ref{ksjhskjsh223} one obtains that the process is ergodic and
weakly mixing.

It follows from theorem 3.4.1 of \cite{PrZa96} that for $\varphi\in
L^2(\mu)$ there exists a set $I\subset [0,+\infty)$ of relative
measure $1$ such that
\begin{equation}
\lim_{|t|\rightarrow \infty,\; t\in
I}\E[\varphi(x_t,\theta_t,v_t,\omega_t)]=\mu[\varphi]\quad
\text{in}\quad L^2(\mu)\text{.}
\end{equation}
Furthermore since $t\rightarrow
\E[\varphi(x_t,\theta_t,v_t,\omega_t)]$ is continuous when $\varphi$
is continuous and bounded we deduce that when $\varphi$ is
continuous and bounded then
\begin{equation}
\lim_{t \rightarrow
\infty}\E[\varphi(x_t,\theta_t,v_t,\omega_t)]=\mu[\varphi]\quad
\text{in}\quad L^2(\mu)\text{.}
\end{equation}
The fact that the process is strongly mixing then follows from
corollary 3.4.3 of \cite{PrZa96}.
\end{proof}

In \cite{BoOw2007c}, using thoerem~\ref{thm:slidingdiskerggodic} we
prove that if $U$ is non-constant then the $x$-displacement of the
sliding disk is $\mu$ a.s.~not ballistic (see
Proposition~\ref{thm:nondirected}).   However, the mean-squared
displacement with respect to the invariant law is ballistic (see
theorem~\ref{thm:ballistic}).  More precisely, we show that the
squared standard deviation of the $x$-displacement with respect to
its noise-average grows like $t^2$. This implies that the process
exhibits not only ballistic transport but also ballistic diffusion.
If $U$ is constant then the squared standard deviation of the
$x$-displacement is diffusive (grows like $t$).   See below for
theoretical results and numerical experiments using efficient
stochastic variational integrators.

\begin{Proposition} \label{thm:nondirected}
Provided that $U$ is non-constant, then $\mu$ a.s.~
\[
\lim_{t \to \infty} \frac{x(t) - x(0)}{t} \to 0 \text{.}
\]
\end{Proposition}

\begin{Proposition}
The squared standard deviation of the $x_t + \theta_t$-degree of
freedom is diffusive, i.e.,
\begin{equation}\label{kskjhs3}
\lim_{t\rightarrow
\infty}\frac{\E_\mu[(x_t+\theta_t-\E[x_t+\theta_t])^2]}{t}=\frac{2
\alpha^2 \sigma^2}{c^2 (\sigma^2+1) } \text{.}
\end{equation}
\end{Proposition}

\begin{Proposition}
Assume that $U$ is non constant, then
\begin{equation}\label{ksddssdkjhs3}
\lim\sup_{t\rightarrow \infty}\frac{\E_\mu\Big[\big(-x_t+\sigma
\theta_t-\E[-x_t+\sigma \theta_t]\big)^2\Big]}{t^2}\leq
4\frac{1+\sigma}{\beta}
\end{equation}
and
\begin{equation}\label{ksddssjkodkjhs3}
\lim\inf_{t\rightarrow \infty}\frac{\E_\mu\Big[\big(-x_t+\sigma
\theta_t-\E[-x_t+\sigma \theta_t]\big)^2\Big]}{t^2}\geq
\frac{1}{4}\frac{1+\sigma}{\beta}
\end{equation}
\end{Proposition}

\begin{Theorem}\label{thm:ballistic}
We have (\cite{BoOw2007c})
\begin{itemize}
\item If $U$ is constant then
\begin{equation}\label{kskjkjhs3}
\lim_{t\rightarrow
\infty}\frac{\E_\mu\big[(x_t-\E[x_t])^2\big]}{t}=\frac{2 \alpha^2
\sigma^2}{c^2 (\sigma^2+1) (\sigma+1)^2}
\end{equation}
\item If $U$ is non constant then
\begin{equation}\label{ksddssdkkjkjhs3}
\lim\sup_{t\rightarrow
\infty}\frac{\E_\mu\big[(x_t-\E[x_t])^2\big]}{t^2}\leq
\frac{4}{\beta(1+\sigma)}
\end{equation}
and
\begin{equation}\label{ksddsdssdkjhs3}
\lim\inf_{t\rightarrow
\infty}\frac{\E_\mu\big[(x_t-\E[x_t])^2\big]}{t^2}\geq
\frac{1}{4\beta(1+\sigma)}
\end{equation}
\end{itemize}
\end{Theorem}

\begin{figure}[htbp]
\begin{center}
\includegraphics[scale=0.3,angle=0]{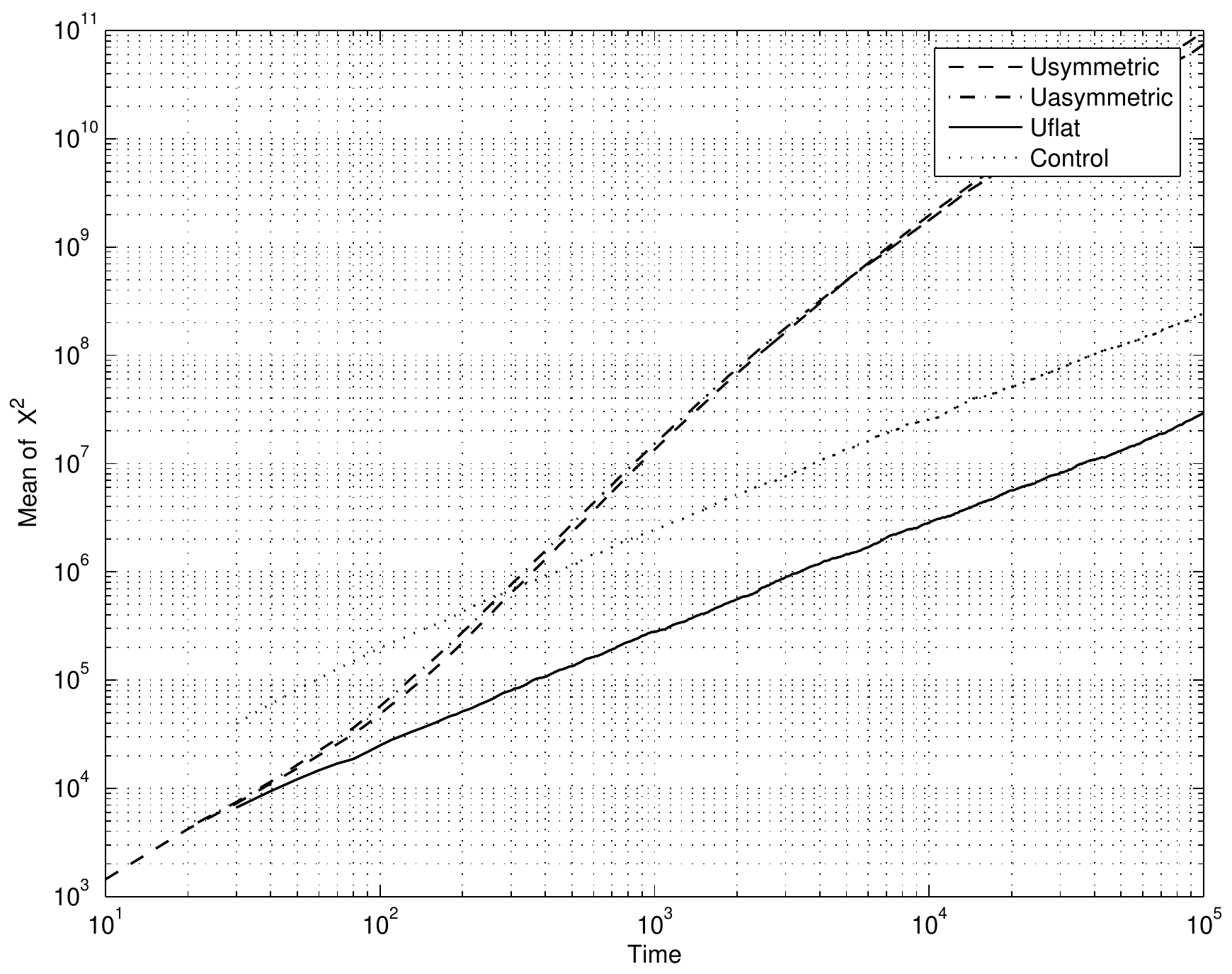}
\caption{\footnotesize {\bf Sliding Disk at Uniform temperature,
$h=0.01$, $\alpha=5.0$, $c=0.1$.} A log-log plot of the mean squared
displacement of the ball.  It clearly shows that the x-position
exhibits anomalous diffusion when $U$ is symmetric or asymmetric.
The disk is started from rest. In the control and flat $U$ cases the
diffusion is normal. } \label{fig:loglogX2}
\end{center}
\end{figure}

\begin{figure}[htbp]
\begin{center}
\includegraphics[scale=0.2,angle=0]{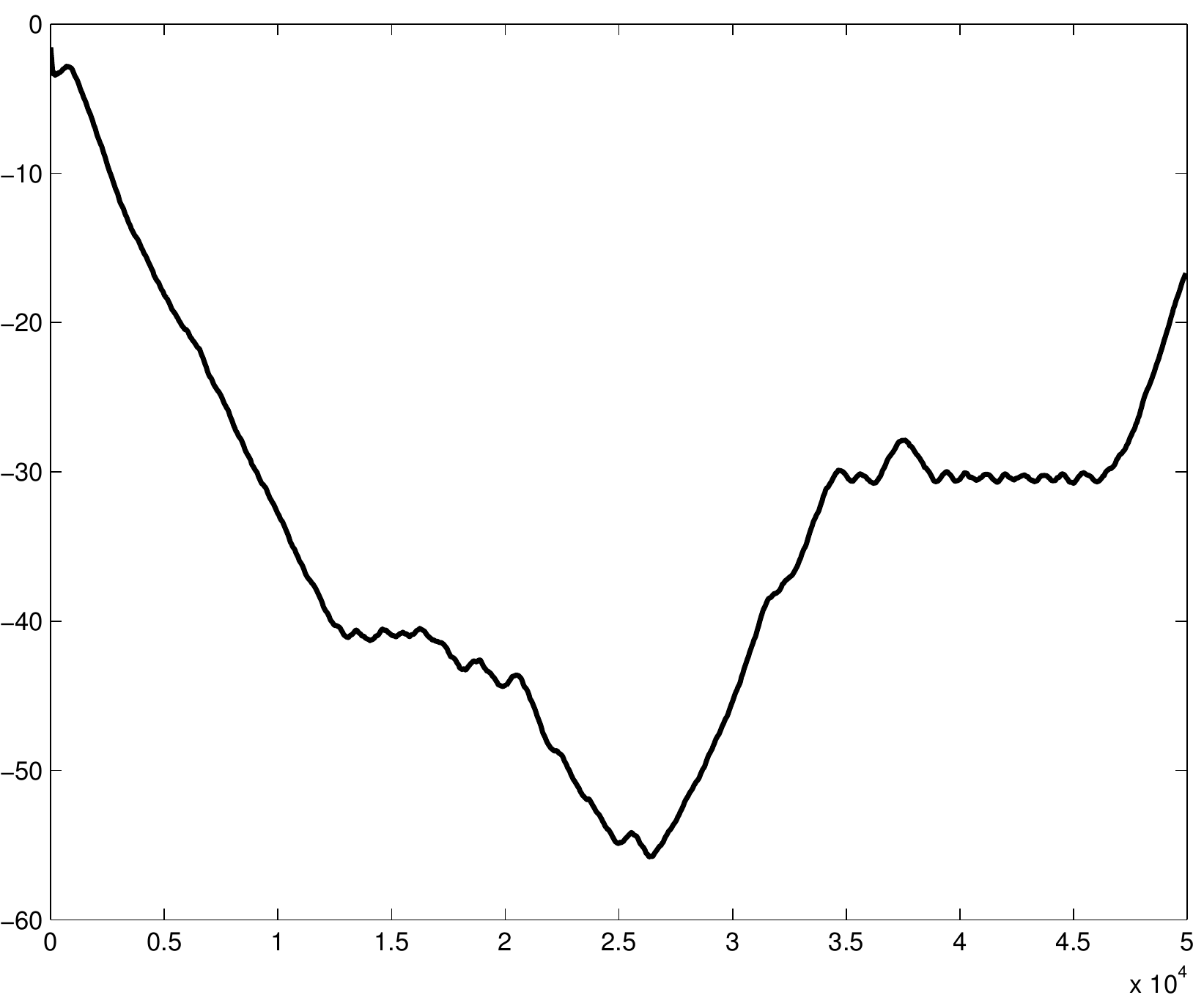}
\includegraphics[scale=0.2,angle=0]{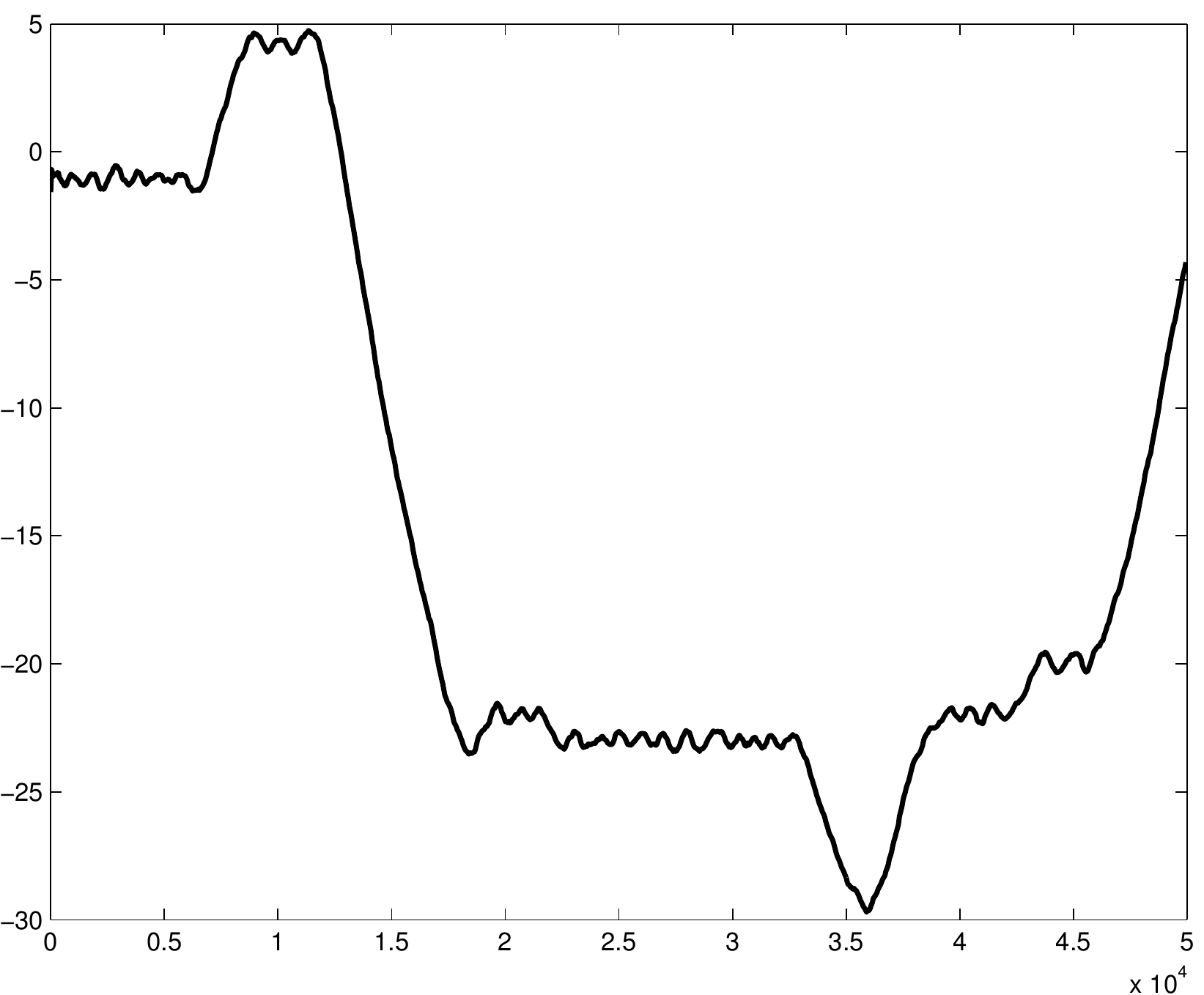}
\includegraphics[scale=0.2,angle=0]{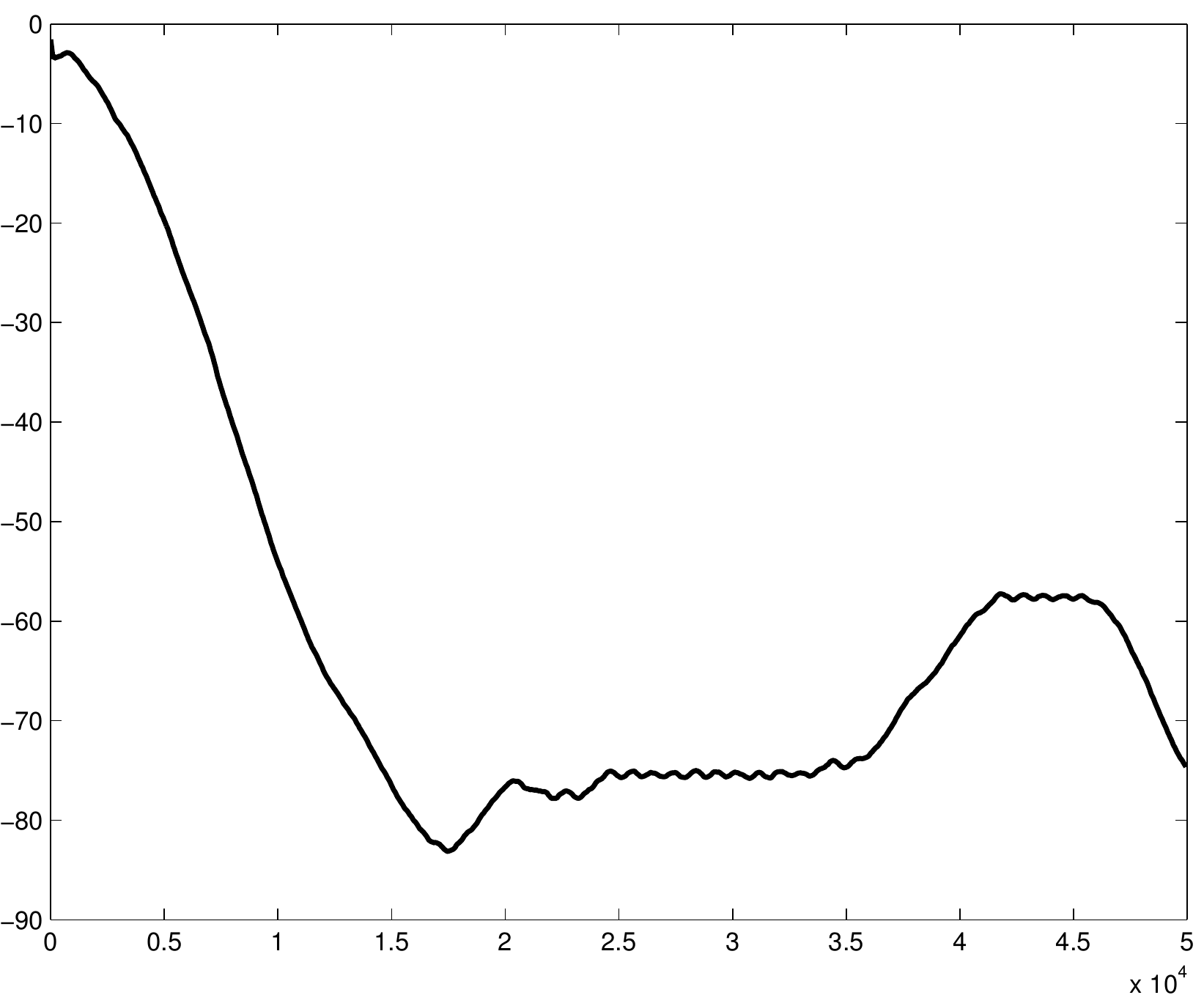}
\includegraphics[scale=0.2,angle=0]{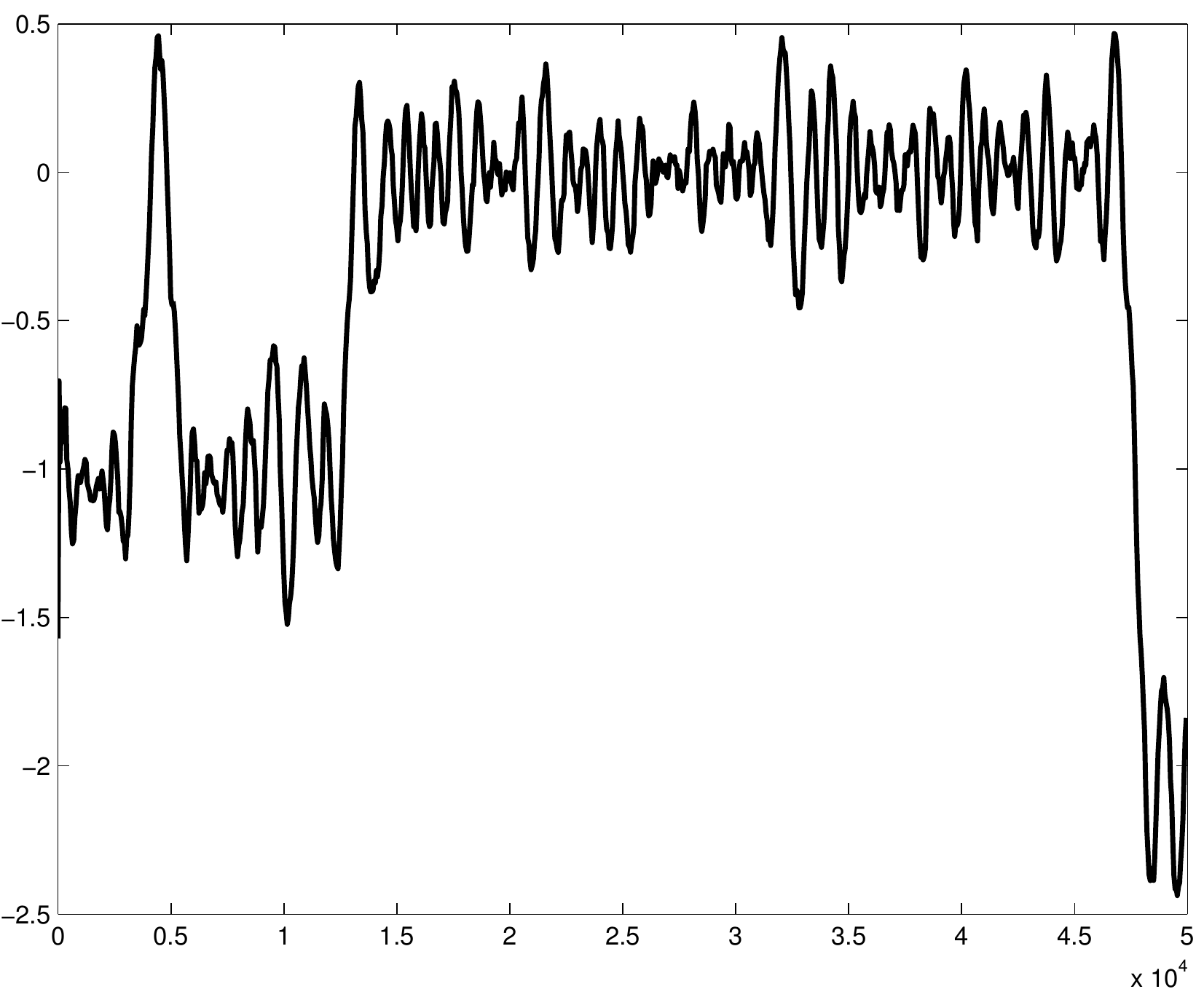} \\
\hbox{\hspace{0.15in}   (a)  $\alpha=0.001$ \hspace{0.3in} (b)
$\alpha=0.001$ \hspace{0.3in} (c)  $\alpha=0.00075$ \hspace{0.3in}
(d) $\alpha=0.00075$  \hspace{0.3in}    } \caption{\footnotesize
{\bf Angular position of a magnetic motor (uniform temperature).}
Four different realizations of the angular component of the center
of mass of a magnetic motor are plotted. The system is started from
rest.
 }
\label{fig:ctmagmotortheta}
\end{center}
\end{figure}

Classical homogenization techniques can't be applied to obtain
theorem \ref{thm:ballistic} (since the stochastic forcing is
degenerate on momentums). We refer to \cite{BoOw2007c} for a proof
of that theorem. The ballistic diffusion is caused by long time
memory effects created by the degeneracy of the noise and the
coupling between the two degrees of freedom through $U$. Figure
\ref{fig:loglogX2} gives an illustration of the mean-squared
displacement of the rolling disk versus time started from rest. In
\cite{BoOw2007c} we have used that phenomenon to propose a
fluctuation driven magnetic motor characterized by ballistic
diffusion at uniform. A plot of the angular displacement of that
magnetic motor versus time for a single realization started from
rest is given in figure \ref{fig:ctmagmotortheta}.

\end{document}